\documentclass[12pt,twoside]{amsart}
\usepackage{amsmath, amsthm, amscd, amsfonts, amssymb, graphicx}
\usepackage[bookmarksnumbered, plainpages]{hyperref}

\newcommand{\pt}{\partial}
\newcommand{\bp}{\overline{\partial}}
\newcommand{\lmd}{\lambda}

\newcommand{\og}{\omega}

\newcommand{\tr}{\text{tr}}

\newcommand{\tend}{\text{End}}

\newtheorem{thm}{Theorem}[section]
\newtheorem{cor}[thm]{Corollary}
\newtheorem{lem}[thm]{Lemma}
\newtheorem{prop}[thm]{Proposition}
\newtheorem{defn}[thm]{Definition}

\numberwithin{equation}{section}

\begin{document}

\title{\bf Analytically stable Higgs bundles on some non-K\"ahler manifolds}
\author{Chuanjing Zhang and Xi Zhang}

\thanks{{\scriptsize
\hskip -0.4 true cm \textit{2010 Mathematics Subject Classification:}
53C07; 14J60; 32Q15.
}}

\maketitle

\begin{abstract}
In this paper, we study Higgs bundles on non-compact Hermitian manifolds. Under some assumptions for the underlying Hermitian manifolds which are not necessarily K\"ahler, we solve the Hermitian-Einstein equation on analytically stable Higgs bundles. 
\end{abstract}

\vskip 0.2 true cm


\pagestyle{myheadings}
\markboth{\rightline {\scriptsize C. Zhang et al.}}
         {\leftline{\scriptsize Higgs bundles on non-K\"ahler manifolds}}

\bigskip
\bigskip


\section{ Introduction}

Let $(M, \omega )$ be an $n$-dimensional  Hermitian manifold and $(E, \overline{\partial}_{E})$ a $r$-rank holomorphic vector bundle on $M$.
A Hermitian metric $H$ on the bundle $E$ is called  $\omega$-Hermitian-Einstein    if it
  satisfies the following Hermitian-Einstein equation on $M$, i.e.
\begin{equation}\label{HE}
\sqrt{-1}\Lambda_{\omega} (F_{H}-\frac{1}{r}\tr F_{H} Id_{E})
=0,
\end{equation}
where  $F_{H}$ is the curvature tensor of Chern connection $D_{H}$ with respect to  $H$ and $\Lambda_{\omega }$ denotes the contraction with  the Hermitian metric $\omega $.

When $(M, \omega)$ is a compact K\"ahler manifold, by the famous Donaldson-Uhlenbeck-Yau theorem (\cite{NS65, DON85, UY86}), we know the holomorphic vector bundle $(E, \overline{\partial}_{E})$ must have an   $\omega$-Hermitian-Einstein metric if $(E, \overline{\partial}_{E})$ is $\omega$-stable  in the sense of Mumford-Takemoto. This classical result has a lot of interesting and important generalizations and extensions (see  \cite{LY, HIT, SIM, br, BS, GP,  HL, ag, Bi, BT, JZ, L1, LN2, M, Mo1, Mo2}, etc.).

  A Higgs bundle $(E , \overline{\partial }_{E}, \theta )$ over $M$ is a holomorphic bundle $(E , \overline{\partial }_{E})$ coupled with a Higgs field $\theta \in \Omega_X^{1,0}(\mathrm{End}(E))$
such that $\overline{\partial}_{E}\theta =0$ and $\theta \wedge \theta =0$.
Higgs bundles  first emerged thirty
years ago in  Hitchin's (\cite{HIT}) reduction of self-dual equation on $R^{4}$ to Riemann surface  and in  Simpson's (\cite{SIM})
 work on nonabelian
Hodge theory, they have  rich structures and
play an important role in many different areas including gauge theory,
K\"ahler and hyperk\"ahler geometry, group representations and
nonabelian Hodge theory.  Letting $H$ be a Hermitian metric on the bundle $E$, we consider the Hitchin-Simpson connection:
 $
 D_{H,  \theta }= D_{H}+\theta + \theta^{*_H}
 $,
where  $\theta^{*_H}$ is the adjoint of $\theta $ with respect to the metric $H$.
The curvature of this connection is
\begin{equation*}
F_{H,\theta}=F_H+[\theta,\theta^{*_H}]+\partial_H\theta+\bar{\partial}_E\theta^{*_H},
\end{equation*}
where  $\partial_H$ is the $(1, 0)$-part of $D_{H}$.  A Hermitian metric $H$  is said to be
a Hermitian-Einstein metric on Higgs bundle $(E, \overline{\partial }_{E}, \theta  )$ if  it satisfies
\begin{equation}\label{HEE}
\sqrt{-1}\Lambda_{\omega} F_{H, \theta }^\bot
=0.
\end{equation}
where $F_{H, \theta }^\bot$ is the trace-free part of the curvature of the Hitchin-Simpson connection.

The Donaldson-Uhlenbeck-Yau theorem was generalized to the Higgs bundles case by Hitchin (\cite{HIT}) and Simpson (\cite{SIM}, \cite{SIM2}). Simpson (\cite{SIM}) even studied some non-compact K\"ahler manifolds case. Under some assumption for the base manifold, he  proved that the analytic stability implies the existence of Hermitian-Einstein metric. The Donaldson-Uhlenbeck-Yau theorem for the non-compact base manifold case is important and interesting (\cite{L1, Mo1, Mo2, Mo4}). Recently, Mochizuki (\cite{Mo3}) made an important progress in this direction. He weakened the assumption in Simpson's result such that the volume of base manifold may not be finite
and he also studied the curvature decay of the Hermitian-Einstein metrics.

In this paper, we  study the  non-K\"ahler case. A Hermitian metric $\omega $ is called to be Gauduchon if it satisfies $\pt\bp \og^{n-1}=0$. If $M$ is compact, it has been proved by Gauduchon (\cite{Gaud}) that  there exists a Gauduchon metric  in the conformal class of every Hermitian metric $\omega$.
When the base Hermitian manifold is compact and Gauduchon, the Donaldson-Uhlenbeck-Yau theorem is  also valid (see \cite{bis0,bis1,Bu, LY, LT,LT95}).
Inspired by Mochizuki's result (\cite{Mo3}), we  consider the case that the base manifold $(M, \omega)$ is non-compact Gauduchon and satisfies the following assumption.

\medskip

{\bf Assumption 1.} Let $\varphi $ be a nonnegative function on $(M, \omega )$ with $\int_{M}\varphi \frac{\omega^{n}}{n!}<+\infty$. There exist positive constants $C_{i}$ (i=1, 2) such that for any nonnegative bounded function $f$ satisfying
\begin{equation}
\sqrt{-1}\Lambda_{\omega}\partial \overline{\partial }f\geq -B\varphi
\end{equation}
in weakly sense (see definition (\ref{def1}) for details) for a positive number $B$, we have
\begin{equation}
\sup_{x\in M} f(x)\leq C_{1}B+C_{2}\int_{M}f\varphi \frac{\omega^{n}}{n!}.
\end{equation}
Moreover, if the function $f$ satisfies
$\sqrt{-1}\Lambda_{\omega}\partial \overline{\partial }f\geq 0$ on $M$, then we have $\sqrt{-1}\Lambda_{\omega}\partial \overline{\partial }f\equiv 0$.

\medskip

Let the background metric $H_0$ be a Hermitian metric of $E$ such that
\begin{equation}\label{A2}|\sqrt{-1}\Lambda_\omega F_{H_0, \theta }|_{H_{0}}\leq \hat{B}\varphi\end{equation} for some constant $\hat{B}>0$.
Define the analytic degree of $E$ to be the real number
\begin{equation*}
\textmd{deg}_{\omega}(E, H_{0})=\sqrt{-1}\int_M\textmd{tr}(\Lambda_{\omega} F_{H_{0},\theta})\frac{\omega^{n}}{n!}.
\end{equation*}
As in \cite{SIM}, we  define the analytic degree of any saturated sub-Higgs sheaf $\mathcal{S}$ of $(E , \overline{\partial }_{E}, \theta )$ by
\begin{equation} \label{cw}
\textmd{deg}_{\omega}(\mathcal{S}, H_{0})=\int_{M\setminus \Sigma_{\mathcal{S}}} (\sqrt{-1}\textmd{tr}(\pi_{\mathcal{S}} \Lambda_{\omega} F_{H_{0},\theta})-|\overline{\partial }_{\theta}\pi_{\mathcal{S}} |_{H_{0}}^{2})\frac{\omega^{n}}{n!},
\end{equation}
where $\Sigma_{\mathcal{S}}$ denotes the set of singularities where $\mathcal{S}$ is not locally free, $\overline{\partial}_{\theta}:=\overline{\partial}_{E}+\theta$ and $\pi_{\mathcal{S}} $ denotes the projection onto $\mathcal{S}$ with respect to the metric $H_{0}$ outside $\Sigma_{\mathcal{S}}$. When the base Gauduchon manifold $(M, \omega)$ is compact, it is easy to see that the analytic degree $\textmd{deg}_{\omega}(\mathcal{S}, H_{0})$ is independent of the choice of the background metric $H_{0}$.


Following \cite{SIM}, we say that the Higgs bundle $(E , \overline{\partial }_{E}, \theta )$  is $H_{0}$-analytic stable (semi-stable) if for every proper saturated sub-Higgs sheaf $\mathcal{S}\subset E$, it holds
\begin{equation}
\frac{\textmd{deg}_{\omega}(\mathcal{S},H_{0})}{\textmd{rank}(\mathcal{S})}<(\leq)\frac{\textmd{deg}_{\omega}(E,H_{0})}{\textmd{rank}(E)}.
\end{equation}
Now we give our main theorem as follows.

\medskip

\begin{thm}\label{theorem1}
Let $(M, \omega )$ be a non-compact Gauduchon manifold satisfying the assumption 1,  $(E , \overline{\partial }_{E}, \theta )$ a Higgs bundle over $M$ and $H_{0}$ a background Hermitian metric on $E$ satisfying the condition (\ref{A2}). If $(E , \overline{\partial }_{E}, \theta )$ is analytically $H_0$-stable, then there exists a Hermitian-Einstein metric $H$ satisfying the below conditions:
\begin{enumerate}
\item  $\det(H)=\det(H_0)$.
\item Set $h=H_0^{-1}H$. Then, $|h|_{H_0}$ and $|h^{-1}|_{H_0}$ are bounded, and $\int_{M} (|  \overline{ \partial}_{E}h|_{H_{0}}^2 +|[\theta , h]|_{H_0}^{2})\frac{\omega^{n}}{n!}< +\infty$.
\end{enumerate}
\end{thm}

The above theorem can be seen as a  generalization of Mochizuki's result (\cite{Mo3}) to the non-K\"ahler case. In \cite{Mo3}, Mochizuki proved the existence of an exhaustion function $\phi$ on $M$.
Fix a number $a_i$ and let $M_i$ denote the compact space $\phi(x)\leq a_i$ with boundary $\partial M_i$, so we can  take a sequence of exhaustion compact subsets $M_i$ in $M$ with $\cup M_i=M$. Let's consider the Dirichlet problem on $M_i$:
\begin{equation}\label{DP}
\left\{\begin{split}
&\sqrt{-1}\Lambda_\omega F_{H_i}^\bot=0,\\
&H_i|_{\partial M_i}= H_0.
\end{split}\right.
\end{equation}
 According to the results of Donaldson (\cite{DON92}, \cite{Z} for the Hermitian manifold case),  we know  that there exists a unique Hermitian metric $H_i$  satisfying the above Dirichlet problem (\ref{DP}) and $\det(H_i)=\det(H_{0})$ on $M_i$. Following the idea in \cite{Mo3}, one can take the limit as $a_i\rightarrow +\infty$ and  get the convergence $H_\infty$ of a subsequence of $H_i$ on any compact subset, which should satisfy  $\sqrt{-1}\Lambda_\omega F_{H_\infty}^\bot=0$ on the whole $M$. Now the key is to obtain a $C^0$-bound. When the base manifold is K\"ahler, Mochizuki (\cite{Mo3}) introduced the Donaldson's functional on the space of Hermitian metrics satisfying the Dirichlet boundary condition. The Donaldson's functional played a key role in Mochizuki's proof of the uniform $C^0$-bound. However, in the non-K\"ahler case, the Donaldson
functional may not be well-defined. So we
need new argument in our case. In fact, our argument relies on the following identity:
\begin{equation}\label{eqn04021}
\int_M \mathrm{tr}(\Phi(H_0,\theta)s)\frac{\omega^n}{n!}+\int_{M}\langle \Psi(s)(\overline{\partial}_{\theta}s),\overline{\partial}_{\theta}s\rangle_{H_0}\frac{\omega^n}{n!}=\int_M \mathrm{tr}(\Phi(H,\theta)s)\frac{\omega^n}{n!},
\end{equation}
where $s=\log (H_{0}^{-1}H)$,
\begin{equation} \Phi(H, \theta)=\sqrt{-1}\Lambda_{\omega } (F_{H}+[\theta,\theta^{*_H}]) \end{equation} and
\begin{equation}\label{eq3301}
\Psi(x,y)=
\left\{\begin{split}
&\frac{e^{y-x}-1}{y-x},\ \ \ &x\neq y;\\
&\ \ \ \  1,\ \ \ \ \ \  &x=y.
\end{split}
\right.
\end{equation}
The above identity (\ref{eqn04021}) was proved in \cite{NZ} for the closed Gauduchon manifold case, and in \cite{ZZZ} for the compact Gauduchon manifold with non-empty boundary and some non-compact case.

This paper is organized as follows. In Section 2,  we give some estimates and
preliminaries for the Hermitian-Einstein equation (\ref{HEE}). In Section 3,  we give a proof of Theorem \ref{theorem1} by using the identity (\ref{eqn04021}). In Section 4,  we study the uniqueness of Hermitian-Einstein metric in Theorem \ref{theorem1}.

\medskip

{\bf  Acknowledgement:} The two authors are partially supported by NSF in China No.11625106, 11571332 and 11721101. The first author is also supported by  NSF in China No.11801535, the China Postdoctoral Science Foundation (No.2018M642515) and the Fundamental Research Funds for the Central Universities.

\medskip

\section{Preliminary results}

Let $(M, \omega )$ be a Hermitian manifold and $(E , \overline{\partial }_{E}, \theta )$ a Higgs bundle over $M$. Letting $H_{0}$ and $H$ be two Hermitian metrics on the bundle $E$, we denote
\begin{equation}\label{7.2}
\mathcal{S}_{H_{0}}(E)=\{\eta \in \Omega^{0}(M, End(E))| \quad \eta ^{\ast H_{0}}=\eta \},
\end{equation}
and set
\begin{equation}
h=H_{0}^{-1}H=\exp{s},
\end{equation}
where $s\in \mathcal{S}_{H_{0}}(E)\cap \mathcal{S}_{H}(E)$. It is easy to check the following identities
\begin{equation}\label{id1}
\begin{split}
&\log (\frac{1}{2r}(\tr h + \tr h^{-1}))\leq |s|_{H_{0}}\leq r^{\frac{1}{2}}\log (\tr h + \tr h^{-1});\\
&\partial _{H}-\partial_{H_{0}} =h^{-1}\partial_{H_{0}}h ;\\
&F_{H}-F_{H_{0}}=\overline{\partial }_{E} (h^{-1}\partial_{H_{0}} h) ;\\
&\theta^{\ast H}=h^{-1} \theta^{\ast H_{0}} h ,\\
\end{split}
\end{equation}
where $r=rank (E)$ and $\partial_{H_{0}}$ is the $(1, 0)$ part of the Chern connection $D_{H_{0}}$. Furthermore,
we have the following estimates (Lemma 3 (d) in \cite{SIM}, Proposition 2.7 in \cite{Z})
\begin{equation}\label{k01}
\sqrt{-1}\Lambda_{\omega}\partial \overline{\partial } \log (\tr h) \geq -|\Phi (H_{0}, \theta)|_{H_{0}}-|\Phi (H, \theta)|_{H}
\end{equation}
and
\begin{equation}\label{k01}
\sqrt{-1}\Lambda_{\omega}\partial \overline{\partial } \log (\tr h + \tr h^{-1}) \geq -|\Phi (H_{0}, \theta)|_{H_{0}}-|\Phi (H, \theta)|_{H}.
\end{equation}

\medskip

The Dirichlet problem for the Hermitian-Einstein equation was first solved in \cite{DON92} by Donaldson for the K\"ahler manifold case,  in \cite{Z} for the general Hermitian manifold case. The following proposition was proved in \cite{ZZZ}.

\medskip

\begin{prop} [Theorem 4.1 in \cite{Z}] \label{comthm}
Let $(E,\bar{\partial}_E,\theta)$ be a Higgs bundle
over a compact Hermitian manifold $(\overline{X}, \omega )$ with non-empty
boundary $\partial X$ and $H_{0}$ a Hermitian metric on $E$. There is a unique
Hermitian metric $H$ on $E$ such that
\begin{equation}\label{DP2}
\left\{\begin{split}
&\sqrt{-1}\Lambda_{\omega } (F_{H}+[\theta,\theta^{*_H}])=\lambda Id_{E},\\
&H |_{\partial X}= H_0,
\end{split}\right.
\end{equation}
where $\lambda $ is a constant.
\end{prop}

\medskip

Let $\tilde{H}$ be a solution of (\ref{DP2}),  $f=\log \frac{\det (H_{0})}{\det (\tilde{H})}$ and $H:=e^{\frac{f}{r}}\tilde{H}$. It is easy to check that
\begin{equation}
\sqrt{-1}\Lambda_{\omega } (F_{H}+[\theta,\theta^{*_H}])=(\lambda -\frac{\sqrt{-1}}{r}\Lambda_{\omega }\partial \overline{\partial } f) Id_{E}.
\end{equation}
So we have:
\begin{equation}\label{DP4}
\left\{\begin{split}
&\sqrt{-1}\Lambda_{\omega } (F_{H}+[\theta,\theta^{*_H}])^\bot=0,\\
&H |_{\partial X}= H_0.
\end{split}\right.
\end{equation} Conversely, if we have  a solution of (\ref{DP4}), then we can get a solution of (\ref{DP2}) by conformal transformation. Using the uniqueness of (\ref{DP2}), one can easily prove  the uniqueness of (\ref{DP4}). So we obtain the following corollary.

\medskip

\begin{cor} \label{DP3}
Let $(E,\bar{\partial}_E,\theta)$ be a Higgs bundle
over a compact Hermitian manifold $(\overline{X}, \omega )$ with non-empty
boundary $\partial X$ and $H_{0}$ a Hermitian metric on $E$. There is a unique
Hermitian metric $H$ on $E$ with $det (H_{0}^{-1}H)=1$ on $\overline{X}$ satisfying
(\ref{DP4}).

\end{cor}

\medskip

 Given $\eta \in \mathcal{S}_{H_{0}}(E)$,  we can choose  a local unitary basis  $\{e_{\alpha}\}_{\alpha=1}^{r}$ with respect to $H_{0}$ and local functions $\{\lmd_{\alpha}\}_{\alpha=1}^{r}$ such that
\begin{eqnarray*}
\eta=\sum_{\alpha=1}^r \lmd_{\alpha} \cdot e_{\alpha}\otimes e^{\alpha}  ,
\end{eqnarray*}
 where  $\{e^{\alpha}\}_{\alpha=1}^{r}$ denotes the dual basis in $E^*$. Let   $\Psi\in C^{\infty}(\mathbb{R}\times \mathbb{R}, \mathbb{R})$ and $A=\sum\limits^r_{\alpha,\beta=1}A_{\beta}^{\alpha} e_{\alpha}\otimes e^{\beta}\in \tend(E)$.  We define:
\begin{equation*}
\Psi(\eta)(A)=\Psi(\lmd_{\beta},\lmd_{\alpha})A^{\alpha}_{\beta} e_{\alpha}\otimes e^{\beta}.
\end{equation*}

\medskip

\begin{prop} [Proposition 2.6 in \cite{ZZZ}] \label{idbundle01}
Let $(E,\bar{\partial}_E,\theta)$ be a Higgs bundle with a fixed Hermitian metric $H_0$
over a compact Gauduchon manifold $(\overline{X}, \omega )$ with non-empty smooth boundary $\partial X$. Let $H$ be a Hermitian metric on $E$ satisfying $H|_{\partial X}=H_0|_{\partial X}$.
Then we have the following identity:
\begin{equation}\label{eq04021}
\int_{\overline{X}} \mathrm{tr}(\Phi(H_0,\theta)s)\frac{\omega^n}{n!}+\int_{\overline{X}}\langle \Psi(s)(\overline{\partial}_{\theta}s),\overline{\partial}_{\theta}s\rangle_{H_0}\frac{\omega^n}{n!}=\int_{\overline{X}} \mathrm{tr}(\Phi(H,\theta)s)\frac{\omega^n}{n!},
\end{equation}
where $s:=\log(H^{-1}_0H)$ and $\Psi $ is the function which is defined in (\ref{eq3301}).
\end{prop}

\medskip

\medskip

In the following, we always assume that the compact manifold $\overline{X}$ with non-empty smooth boundary $\partial X$ is a subset of the Hermitian manifold $(M, \omega )$.  Let $H$ be the unique Hermitian metric in Corollary \ref{DP3}. Set $\exp{s}=h=H_{0}^{-1}H$. By the condition $\det (h)\equiv 1$ on $\overline{X}$ and the relationship between the geometric mean and arithmetic mean, one can get that
 \begin{equation}\label{K02}
\frac{1}{r}\tr(\exp s)\geq (\det h)^{\frac{1}{r}}=1
\end{equation}
on $\overline{X}$. As that in  \cite{Mo3}, we
  extend $\log(\tr(\exp s)/r)$ and $(|\Lambda_\omega F_{H_{0}, \theta }|_{H_{0}})|_{\overline{X}}$ to the functions $\log(\tr(\exp s)/r)^{\thicksim}$ and $|\Lambda_\omega F_{H_{0}, \theta }|_{H_{0}}^{\thicksim X}$ on the whole $M$ by setting $0$ outside $\overline{X}$.

Let $g$ be the Hermitian metric with respect to $\omega $. As usual, we denote the Beltrami-Laplace operator on the Hermitian manifold $(M, \omega )$ by $\Delta_{g }$, and
define the complex Laplace operator $\widetilde{\Delta}_{\omega }$ for functions as
\[
\widetilde{\Delta}_{\omega }f=2\sqrt{-1}\Lambda_{\omega } \partial \overline{\partial }
f .
\]
 It is well known that the difference of the two Laplacians is given by a first order differential operator as follows
\begin{equation*} \label{laplacian}
(\widetilde{\Delta}_{\omega }-\Delta_{g})f=\langle V,\nabla f\rangle_g,
\end{equation*}
where $V$ is a smooth vector field on $M$. Usually the complex Laplace operator is not a self adjoint operator.

\medskip

\begin{defn}\label{def1} A function $f$ on the Hermitian manifold $(M, \omega )$ satisfying
\begin{equation}
\sqrt{-1}\Lambda_{\omega}\partial \overline{\partial }f\geq \eta
\end{equation}
in weakly sense means that, for any nonnegative compactly supported smooth function $\psi$, there holds
\begin{equation}
\int_{M}f\sqrt{-1}\Lambda_{\omega}\partial \overline{\partial }(\psi \frac{\omega^{n-1}}{(n-1)!} )\geq \int_{M}\eta \psi \frac{\omega^{n}}{n!}.
\end{equation}
\end{defn}

\medskip

\begin{prop}\label{prop2}
We have
\begin{equation}
\sqrt{-1}\Lambda_\omega\partial\overline{\partial}\log(\tr(\exp s)/r)^{\thicksim}\geq -|\Lambda_\omega F_{H_{0}, \theta}^\bot|_{H_{0}}^{\thicksim X}
\end{equation}
in weak sense on $M$.
\end{prop}

{\bf Proof} Due to the inequality (\ref{K02}) and the boundary condition of $H$, we know
\begin{equation}
\frac{\partial\tr(\exp s) /r}{\partial\nu}\leq 0,
\end{equation}
where  $\nu$ is the outer normal vector field at $\partial{X}$.
Direct computations give us that
\begin{equation}
\begin{split}
&\int_M \log(\tr(\exp s)/r)^{\thicksim}\cdot \sqrt{-1}\partial\overline{\partial}\Big(\psi\frac{\omega^{n-1}}{(n-1)!}\Big)\\
=&\int_{\overline{X}}\log(\tr(\exp s)/r)\cdot \sqrt{-1}\partial\overline{\partial}\Big(\psi\frac{\omega^{n-1}}{(n-1)!}\Big)\\
=&\int_{\overline{X}}\sqrt{-1}\partial\Big(\log(\tr(\exp s)/r)\cdot\overline{\partial}\Big(\psi\frac{\omega^{n-1}}{(n-1)!}\Big)\Big)
-\int_{\overline{X}}\sqrt{-1}\partial\log(\tr(\exp s)/r)\wedge\overline{\partial}\Big(\psi\frac{\omega^{n-1}}{(n-1)!}\Big)\\
=&\int_{\overline{X}}\sqrt{-1}\overline{\partial}\Big(\partial\log(\tr(\exp s)/r)\wedge\psi\frac{\omega^{n-1}}{(n-1)!}\Big)
-\int_{\overline{X}}\sqrt{-1}\overline{\partial}\partial\log(\tr(\exp s)/r)\wedge\psi\frac{\omega^{n-1}}{(n-1)!},
\end{split}
\end{equation}
where $\psi$ is a test function.
This means that
\begin{equation}
\begin{split}
&\int_{\overline{X}}\sqrt{-1}\overline{\partial}\Big(\partial\log(\tr(\exp s)/r)\wedge\psi\frac{\omega^{n-1}}{(n-1)!}\Big)\\
=&\int_{\overline{X}}\log(\tr(\exp s)/r)\cdot \sqrt{-1}\partial\overline{\partial}\Big(\psi\frac{\omega^{n-1}}{(n-1)!}\Big)
-\int_{\overline{X}}\sqrt{-1}\partial\overline{\partial}\log(\tr(\exp s)/r)\wedge\psi\frac{\omega^{n-1}}{(n-1)!}.
\end{split}
\end{equation}
Similarly, we have
\begin{equation}
\begin{split}
&\int_{\overline{X}}\sqrt{-1}\partial\Big(\overline{\partial}\log(\tr(\exp s)/r)\wedge\psi\frac{\omega^{n-1}}{(n-1)!}\Big)\\
=&\int_{\overline{X}}\log(\tr(\exp s)/r)\cdot \sqrt{-1}\overline{\partial}\partial\Big(\psi\frac{\omega^{n-1}}{(n-1)!}\Big)
-\int_{\overline{X}}\sqrt{-1}\overline{\partial}\partial\log(\tr(\exp s)/r)\wedge\psi\frac{\omega^{n-1}}{(n-1)!}.
\end{split}
\end{equation}
Combining the above two equalities, we deduce that
\begin{equation}
\begin{split}
&2\Big(\int_{\overline{X}}\log(\tr(\exp s)/r)\cdot \sqrt{-1}\partial\overline{\partial}\Big(\psi\frac{\omega^{n-1}}{(n-1)!}\Big)
-\int_{\overline{X}}\sqrt{-1}\partial\overline{\partial}\log(\tr(\exp s)/r)\wedge\psi\frac{\omega^{n-1}}{(n-1)!}\Big)\\
=&\int_{\overline{X}}\sqrt{-1}\overline{\partial}\Big(\partial\log(\tr(\exp s)/r)\wedge\psi\frac{\omega^{n-1}}{(n-1)!}\Big)-\int_{\overline{X}}\sqrt{-1}\partial\Big(\overline{\partial}\log(\tr(\exp s)/r)\wedge\psi\frac{\omega^{n-1}}{(n-1)!}\Big)\\
=&\int_{\overline{X}}d\Big(\sqrt{-1}(\partial\log(\tr(\exp s)/r)-\overline{\partial}\log(\tr(\exp s)/r))\wedge\psi\frac{\omega^{n-1}}{(n-1)!}\Big).
\end{split}
\end{equation}

\medskip

Now we should use the following result.

\medskip

\begin{lem}
Let $(M, \omega )$ be a Hermitian manifold, $g$ be the Riemannian metric with respect to $\omega$. Then we have
\begin{equation}\label{key01}
-\ast df=\sqrt{-1}(\partial f-\overline{\partial}f)\wedge\frac{\omega^{n-1}}{(n-1)!},
\end{equation}
where $\ast$ is the Hodge star operator with respect to the Riemannian metric $g$ and $f$ is a differential function on $M$.
\end{lem}

\medskip

{\bf Proof} Let $\theta$ be a $1$-form. By a simple calculation, we obtain that
\begin{equation}
\theta\wedge\ast df=\langle\theta, df\rangle_g\frac{\omega^n}{n!}=\langle\theta^{1,0}, \overline{\partial}f\rangle_g\frac{\omega^n}{n!}+\langle\theta^{0,1}, \partial f\rangle_g\frac{\omega^n}{n!}
\end{equation}
and
\begin{equation}
\begin{split}
&\theta\wedge\sqrt{-1}(\partial f-\overline{\partial}f)\wedge\frac{\omega^{n-1}}{(n-1)!}\\
=&-\theta^{1,0}\wedge\sqrt{-1}\overline{\partial}f\wedge\frac{\omega^{n-1}}{(n-1)!}
+\theta^{0,1}\wedge\sqrt{-1}\partial f\wedge\frac{\omega^{n-1}}{(n-1)!}\\
=&-\langle\theta^{1,0}, \overline{\partial}f\rangle_g\frac{\omega^n}{n!}-\langle\theta^{0,1}, \partial f\rangle_g\frac{\omega^n}{n!}.
\end{split}
\end{equation}
Therefore, because of the arbitrary of $\theta$, the above lemma follows.

Note that
\begin{equation}
div(\nabla f)\frac{\omega^n}{n!}=d\ast(df).
\end{equation}
One can easily check that
\begin{equation}
div(\psi\nabla f)\frac{\omega^n}{n!}=d\ast(\psi df).
\end{equation}
Hence we derive
\begin{equation}
\begin{split}
&2\Big(\int_{\overline{X}}\log(\tr(\exp s)/r)\cdot \sqrt{-1}\partial\overline{\partial}\Big(\psi\frac{\omega^{n-1}}{(n-1)!}\Big)
-\int_{\overline{X}}\sqrt{-1}\partial\overline{\partial}\log(\tr(\exp s)/r)\wedge\psi\frac{\omega^{n-1}}{(n-1)!}\Big)\\
=&-\int_{\overline{X}}d(\ast\psi d\log(\tr(\exp s)/r))\\
=&-\int_{\overline{X}}div(\psi\nabla\log(\tr(\exp s)/r))\frac{\omega^n}{n!}\\
=&-\int_{\partial X}\psi\frac{\partial\log(\tr(\exp s)/r)}{\partial \nu}\frac{\omega^n}{n!}\\
\geq& 0.
\end{split}
\end{equation}
This implies
\begin{equation}
\begin{split}
\int_M \log(\tr(\exp s)/r)^{\thicksim}\cdot \sqrt{-1}\partial\overline{\partial}\Big(\psi\frac{\omega^{n-1}}{(n-1)!}\Big)\geq &\int_{\overline{X}}\sqrt{-1}\partial\overline{\partial}\log(\tr(\exp s)/r)\wedge\psi\frac{\omega^{n-1}}{(n-1)!}\\
\geq &-\int_{\overline{X}}|\Lambda_\omega F_{H_{0}, \theta}^\bot|_{H_{0}}\psi\frac{\omega^n}{n!}\\
= &-\int_M |\Lambda_\omega F_{H_{0}, \theta}^\bot|_{H_{0}}^{\thicksim X}\cdot \psi\frac{\omega^n}{n!}.
\end{split}
\end{equation}
This conclude the proof of the proposition.

\hfill $\Box$ \\

\section{A proof of Theorem \ref{theorem1}}

Let's first recall  the following lemmas which are proved by Mochizuki in \cite{Mo3}.

\medskip

\begin{lem}(Lemma 2.20 in \cite{Mo3})
There exists an exhaustion function $\phi \in C^{\infty}(M)$.
\end{lem}
\begin{lem}(Lemma 2.21 in \cite{Mo3})
There is a sequence of compact subsets $M_i\subset M_{i+1}\subset M(i=1, 2, \ldots)$ with $\cup M_i=M$, where each $M_i$ is a submanifold with non-empty smooth boundary $\partial M_i$ such that $M_i\setminus \partial M_i$ is an open subset of $M$. Moreover, each connected component of $M_i$ has non-empty boundary.
\end{lem}

Fix a number $a_i$ and let $M_i$ denote the compact space $\{x\in M |\phi(x)\leq a_i \}$ with boundary $\partial M_i$. By choosing a sequence $a_{i}\rightarrow +\infty$ such that each $a_{i}$ is not a critical value of $\phi$, we have a sequence of exhaustion  compact subsets $M_i\subset M$ with non-empty smooth boundary. Let's consider the Dirichlet problem on $M_i$. Then corollary \ref{DP3} tells us that there exists a Hermitian metric $H_{i}$ on $E|_{M_{i}}$ such that
\begin{equation}\label{DP6}
\left\{\begin{split}
&\sqrt{-1}\Lambda_\omega F_{H_i, \theta }^\bot=0,\\
&H_i|_{\partial M_i}= H_0,\\
& \det (H_{0}^{-1}H_{i})=1.
\end{split}\right.
\end{equation}

Set $H_0^{-1}H_i=h_i=\exp s_i$. We already know on $M_i$, it holds that
\begin{equation}
\sqrt{-1}\Lambda_\omega\partial\overline{\partial}\log(\tr(\exp s_i)/r)\geq -|\Lambda_\omega F_{H_{0, i}}^\bot|_{H_{0, i}}.
\end{equation}
 Now extend $\log(\tr(\exp s_i)/r)$ and $(|\Lambda_\omega F_{H_{0}, \theta }^\bot|_{H_{0}})|_{M_{i}}$ to the functions $\log(\tr(\exp s_i)/r)^{\thicksim}$ and $|\Lambda_\omega F_{H_{0}, \theta }^\bot|_{H_{0}}^{\thicksim M_{i}}$ on $M$ by setting $0$ outside $M_i$. According to Proposition \ref{prop2}, we obtain
\begin{equation}
\sqrt{-1}\Lambda_\omega\partial\overline{\partial}\log(\tr(\exp s_i)/r)^{\thicksim}\geq -|\Lambda_\omega F_{H_{0}, \theta}^\bot|_{H_{0}}^{\thicksim M_i}
\end{equation}
in weak sense on $M$.
From the assumptions in Theorem \ref{theorem1}, it can be seen that there exists two positive constants $C_1$, $C_2$ such that for any $i$, we have
\begin{equation}
\sup_{M_i}\log(\tr(\exp s_i)/r)\leq C_1 \hat{B}+C_2\int_{M_i} \log(\tr(\exp s_i)/r)\cdot\varphi\frac{\omega^n}{n!}.
\end{equation}
Note that $\tr s_i =0$, it is easy to check that
\begin{equation}
\log(\tr(\exp s_i)/r)\leq |s_i|_{H_0}\leq (r-1)r^{\tfrac{1}{2}}\log(\tr(\exp s_i)).
\end{equation}
So clearly it implies that
\begin{equation}\label{L1}
\sup_{M_i}|s_i|_{H_0}\leq C_3+C_4\int_{M_i}|s_i|_{H_0}\cdot\varphi\frac{\omega^n}{n!},
\end{equation}
where $C_3$ and $C_4$ are positive constants depending only on $C_{1}$, $C_2$, $\hat{B}$ and $r$.

\medskip

{\bf Proof of Theorem \ref{theorem1}}

When the Higgs bundle $(E,\bar{\partial}_E,\theta)$ is $H_0$-analytically stable on $(M, \omega)$, we will show that, by choosing a subsequence, $H_i$ converge to a Hermitian-Einstein metric $H_{\infty}$ in $C_{loc }^{\infty}$-topology as $i \rightarrow +\infty $.

(1) {\bf Uniform $C^{0}$-estimate}.
By (\ref{L1}), the key is to get a uniform estimate for $\int_{M_i}|s_i|_{H_0}\cdot\varphi\frac{\omega^n}{n!}$, i.e. there exists a constant $\hat{C}$ independent of $i$, such that
\begin{equation}\label{L102}
l_{i}:= \int_{M_i}|s_i|_{H_0}\cdot\varphi\frac{\omega^n}{n!}\leq \hat{C}
\end{equation}
for all $i$.

As that in \cite{Mo3}, we prove (\ref{L102}) by contradiction. If not, there would exist a subsequence $i\rightarrow +\infty$ such that
$l_i\rightarrow +\infty$. Set
\begin{equation} u_i=\frac{s_i}{l_i}.\end{equation}
Then, we have
\begin{equation}\label{id11}
\int_{M_i}|u_i|_{H_0}\cdot\varphi\frac{\omega^n}{n!}=1,
\end{equation}
and
\begin{equation} \label{uc0}
\sup\limits_{M_i}|u_i|_{H_0}\leq \frac{1}{l_i}(C_3+C_4l_i)<C_5<+\infty.
\end{equation}

Now, we show that $\| u_i\|_{L^2_1}$ are uniformly bounded on any compact subset of $ M$.

Based on Proposition \ref{idbundle01},  we deduce
\begin{equation} \label{seq1}
\int_{M_i} \textmd{tr}(\sqrt{-1}\Lambda_\omega F_{H_{0}, \theta}^\bot u_i)\frac{\omega^n}{n!}+l_i\int_{M_i}\langle \Psi(l_iu_i)(\overline{\partial}_{\theta}u_i),\overline{\partial}_{\theta}u_i\rangle_{H_{0}}\frac{\omega^n}{n!}=0.
\end{equation}
By the definition (\ref{eq3301}), it is easy to check that
\begin{equation} \label{seq2}
l\Psi(lx,ly)\rightarrow
\left\{\begin{split}
&(x-y)^{-1},\ \ \ &x>y;\\
\ \ &+\infty,\ \ \ &x\leq y,
\end{split}
\right.
\end{equation}
increases monotonically as $l\rightarrow +\infty$.  Let $\varsigma\in C^{\infty} (\mathbb{R}\times \mathbb{R}, \mathbb{R}^+)$ satisfying $\varsigma(x,y)<(x-y)^{-1}$ whenever $x>y$. Clearly Eqs. \eqref{seq1}, \eqref{seq2}  and the arguments in \cite[Lemma 5.4]{SIM} yield that
\begin{equation} \label{seq3}
\int_{M_i} \textmd{tr}(\sqrt{-1}\Lambda_\omega F_{H_{0}, \theta}^\bot u_i)\frac{\omega^n}{n!}+\int_{M_i}\langle \varsigma(u_i)(\overline{\partial}_{\theta}u_i),\overline{\partial}_{\theta}u_i\rangle_{H_0}\frac{\omega^n}{n!}\leq 0, \ \ i\gg 0.
\end{equation}
From  \eqref{uc0}, we may assume that $(x,y)\in (-C_{5},C_{5})\times (-C_{5},C_{5})$. Note that $\frac{1}{2C_5}<\frac{1}{x-y}$ when $x>y$. In particular,  taking $\zeta(x,y)=\frac{1}{2C_5}$ in (\ref{seq3}), we immediately get
\begin{equation}
\int_{M_i}\textmd{tr}(\sqrt{-1}\Lambda_\omega F_{H_{0}, \theta}^\bot u_i)\frac{\omega^n}{n!}
+\frac{1}{2C_5}\int_{M_i}|\overline{\partial}_{\theta}(u_i)|^2_{H_0}\frac{\omega^n}{n!}\leq 0,
\end{equation}
for $i\gg 0$, and then
\begin{equation}
\int_{M_i}|\overline{\partial}_{\theta}(u_i)|^2_{H_0}\frac{\omega^n}{n!}\leq 2C^2_{5}\int_{M_i}|\sqrt{-1}\Lambda_\omega F_{H_{0}, \theta}^\bot |_{H_0}\frac{\omega^n}{n!}\leq C_6,
\end{equation}
where $C_6$ is a uniform constant. Thus, $u_i$ are bounded in $L_1^2$ on any compact subset $M$. By choosing a subsequence, we have  $u_i\rightharpoonup u_{\infty}$ weakly in $L^2_{1, loc}$. Of course $\tr s_i=0$ and (\ref{uc0}) imply that
\begin{equation} \label{uc01}
\tr u_{\infty} =0, \quad \sup\limits_{M}|u_{\infty}|_{H_0}\leq C_5<+\infty.
\end{equation}

The condition $\int_{M}\varphi \frac{\omega^{n}}{n!}<+\infty$ means that, for any $\epsilon >0$, there exists $
i_{0}$ such that
\begin{equation}\label{si01}
0\leq \int_{M\setminus M_i}\varphi \frac{\omega^{n}}{n!}<\epsilon
\end{equation}
for all $i\geq i_{0}$.
 Combining this and (\ref{id11}), (\ref{uc0}) gives that
\begin{equation}
1\geq \int_{M_i} |u_j|_{H_0}\varphi \frac{\omega^n}{n!} =(\int_{M_j}-\int_{M_j\setminus M_i})|u_j|_{H_0}\varphi \frac{\omega^n}{n!}\geq 1-C_5 \epsilon
\end{equation}
for any $i_{0}\leq i \leq j$.
Noting that $L_1^2\hookrightarrow L^1$ on any compact subset, we have
\begin{equation}
1\geq \int_{M_i} |u_{\infty}|_{H_0}\varphi \frac{\omega^n}{n!} \geq 1-C_5 \epsilon
\end{equation}
for all $i_{0}\leq i$, and then
\begin{equation}
1\geq \int_{M} |u_{\infty}|_{H_0}\varphi \frac{\omega^n}{n!} \geq 1-C_5 \epsilon .
\end{equation}
 This indicates that \begin{equation}
 \int_{M} |u_{\infty}|_{H_0}\varphi \frac{\omega^n}{n!} =1
\end{equation} and $u_{\infty}$ is non-trivial.
If $i_{0}\leq i \leq j$, we derive
\begin{equation}
\begin{split}
&\int_{M_i} \textmd{tr}(\sqrt{-1}\Lambda_\omega F_{H_{0}, \theta}^\bot u_j)\frac{\omega^n}{n!}+\int_{M_i}\langle \varsigma(u_j)(\overline{\partial}_{\theta}u_j),\overline{\partial}_{\theta}u_j\rangle_{H_0}\frac{\omega^n}{n!}\\
\leq& (\int_{M_j}-\int_{M_j\setminus M_i})\textmd{tr}(\sqrt{-1}\Lambda_\omega F_{H_{0}, \theta}^\bot u_j)\frac{\omega^n}{n!} + \int_{M_j}\langle \varsigma(u_j)(\overline{\partial}_{\theta}u_j),\overline{\partial}_{\theta}u_j\rangle_{H_0}\frac{\omega^n}{n!}\\
\leq & 2C_5 \hat{B} \epsilon ,\\
\end{split}
\end{equation}
where we have used (\ref{si01}), (\ref{seq3}) and (\ref{A2}). Taking limits $j\rightarrow \infty $ and $i\rightarrow \infty$, one can obtain
\begin{equation}
\int_{M} \textmd{tr}(\sqrt{-1}\Lambda_\omega F_{H_{0}, \theta}^\bot u_{\infty})\frac{\omega^n}{n!}+\int_{M}\langle \varsigma(u_{\infty})(\overline{\partial}_{\theta}u_{\infty}),\overline{\partial}_{\theta}u_{\infty}\rangle_{H_0}\frac{\omega^n}{n!}\leq  2C_5 \hat{B} \epsilon .
\end{equation}
The fact that $\epsilon $  is arbitrary in the above inequality obviously implies
\begin{equation} \label{seq4}
\int_M \textmd{tr}(\sqrt{-1}\Lambda_\omega F_{H_{0}, \theta}^\bot u_{\infty})\frac{\omega^n}{n!}+\int_{M}\langle \varsigma(u_{\infty})(\overline{\partial}_{\theta}u_{\infty}),\overline{\partial}_{\theta}u_{\infty}\rangle_{H_0}\frac{\omega^n}{n!}\leq 0.
\end{equation}

Now following Simpson's argument  \cite[Lemma 5.5]{SIM}, we conclude that the eigenvalues of $u_{\infty}$ are constant almost everywhere. Let $\mu_1<\mu_2<\cdots<\mu_l$ be the distinct eigenvalues of $u_{\infty}$. Because $\textmd{tr}(u_{\infty})=0$ and $u_{\infty}\neq 0$, there must hold that $2\leq l\leq r$. For each $\mu_{\alpha} (1\leq \alpha\leq l-1)$, we construct a function $P_{\alpha}: \mathbb{R}\rightarrow \mathbb{R}$ such that
\begin{equation*}
P_{\alpha}=
\left\{\begin{split}
&1,\ \ \ x\leq \mu_{\alpha};\\
&0,\ \ \ x\geq \mu_{\alpha+1}.
\end{split}
\right.
\end{equation*}
Setting $\pi_{\alpha}=P_{\alpha}(u_{\infty})$, from \cite[p.887]{SIM}, we have: (i) $\pi_{\alpha}\in L^2_1$; (ii)$\pi_{\alpha}^2=\pi_{\alpha}=\pi_{\alpha}^{*_{H_0}}$; (iii) $(\textmd{Id}_E-\pi_{\alpha})\bar{\partial}\pi_{\alpha}=0$ and (iv) $(\textmd{Id}_E-\pi_{\alpha})[\theta, \pi_{\alpha}]=0$. By Uhlenbeck and Yau's regularity statement of $L^2_1$-subbundle \cite{UY86}, $\{\pi_{\alpha}\}_{\alpha=1}^{l-1}$ determine $l-1$ Higgs sub-sheaves $\{\mathcal{E}_{\alpha }\}_{\alpha=1}^{l-1}$ of $E$.  Using the same argument as that (\cite[Proposition 5.3]{SIM}), we can prove that there must exist a Higgs subsheaf $\mathcal{E}_{\alpha }$ which contradicts the stability of $(E,\bar{\partial}_E,\theta)$. This completes the proof of uniform $C^{0}$-estimate.

(2) {\bf Uniform local $C^{1}$-estimate}.
From the property that $H_i$ satisfies the Hermitian-Einstein equation (\ref{HEE}) and $\det{h_i} =1$ on $M_i$, it is easy to see that
\begin{equation}
\tr{F_{H_i}}=\tr{F_{H_{0}}},
\end{equation}
\begin{equation}\label{c1002}
\begin{split}
\sqrt{-1}\Lambda_\omega\overline{\partial}_{E}\partial_{H_{0}}h_i &=-\sqrt{-1}\Lambda_\omega h_i \cdot F_{H_{0}, \theta }^\bot +\sqrt{-1}\Lambda_\omega \overline{\partial}_{E}h_{i} \cdot h_i^{-1} \cdot \partial_{H_{0}}h_i \\
& + \sqrt{-1}\Lambda_\omega \{[\theta , h_i]\cdot h_i^{-1}\cdot [\theta^{\ast H_0} , h_i]-[[\theta^{\ast H_{0}}, h_i], \theta ]\},\\
\end{split}
\end{equation}
and then
\begin{equation}\label{c1003}
\sqrt{-1}\Lambda_\omega \overline{\partial}\partial \tr{h_i}=-\sqrt{-1}\Lambda_\omega \tr({h_i \cdot F_{H_{0}, \theta }^\bot }) -|  h_i^{-\frac{1}{2}}\cdot \partial_{H_0}h_i|_{H_{0}}^2 -|[\theta , h_i]\cdot h_i^{-\frac{1}{2}}|_{H_0}^{2}.
\end{equation}

Let $T_i=h_i^{-1}\partial_{H_0}h_i$. A direct computation gives us that
\begin{equation}\label{c1001}
\begin{split}
&\sqrt{-1}\Lambda_\omega\partial\overline{\partial}|T_i|^2_{H_i} \geq \frac{1}{2}|\nabla_{H_i}T_i|^2_{H_i} -C_7(|F_{H_0}|_{H_i}+|\theta|_{H_i}^{2}+|\sqrt{-1}\Lambda_\omega F_{H_0}|_{H_i}
+|Rm(g)|_{g}+|\nabla_{g}J|)|T_i|^2_{H_i}\\
& -C_8|D_{H_0}(\Lambda_{\omega }F_{H_0})|_{H_i}|T_i|_{H_i }-C_9|\nabla_{H_0}\theta |_{H_i}^{2},\\
\end{split}
\end{equation}
where the constants $C_7, C_8, C_9$ depend only on the dimension $n$ and the rank $r$. We will follow the argument in \protect {\cite[Lemma 2.4]{LZZ}} to get local uniform $C^{1}$-estimate. Let $\Omega$ be a compact subset in $M$, $d$ be a constant less than the distance of $\Omega $ to $\partial M_{i_{0}}$, where $i_{0}$ is large enough such that $\Omega \subset M_{i_{0}}$. Set $\Omega_1=\{x \in M | \mathrm{dist} (x, \Omega )\leq \frac{1}{4}d\}$ and $\Omega_2=\{x \in M | \mathrm{dist} (x, \Omega )\leq \frac{1}{2}d\}$. Let's choose two non-negative cut-off functions  $\psi_1$, $\psi_2$ such that:
 \begin{equation*}
     \psi_1
     =\left\{\begin{split}
       &0,\ \ \  x\in M\backslash \Omega_1, \\
       &1,\ \ \  x\in \Omega,
     \end{split}
     \right.
  \end{equation*}
   \begin{equation*}
    \psi_2
     =\left\{\begin{split}
       &0, \ \ \ x\in M\backslash \Omega_2 , \\
       &1, \ \ \ x\in \Omega_1,
     \end{split}
     \right.
  \end{equation*}
  and
  $$|\textmd{d} \psi_{\alpha } |^2+|\Lambda_\omega\partial\overline{\partial} \psi_{\alpha}|\leq C_{10}, \ \ \alpha=1,2,$$
where $C_{10}$ is a constant depending only on $d^{-2}$ and the geometry of $(\Omega_2, \omega)$.
Consider the following test function
\begin{equation}\eta_{i} =\psi_1^2|T_i|^2_{H_i}+\tilde{B}\psi_2^2\textmd{tr}h_i,\end{equation}
where the constant $\tilde{B}$ will be chosen large enough later and $i\geq i_{0}$.

It follows from \eqref{c1003} and \eqref{c1001} that
\begin{equation}\sqrt{-1}\Lambda_\omega\partial\overline{\partial}\eta_i\geq \psi^2_2(\tilde{B}e^{-C_{11}}-C_{12})|T_i|^2_{H_i}-C_{13},\end{equation}
where $C_{11}$ is a positive constant depending only on $\sup_{\Omega_2}(\tr{h_i}+\tr{h_i^{-1}})$, $C_{12}$ and $C_{13}$ are positive constants depending only on $\sup_{\Omega_2}(\tr{h_i}+\tr{h_i^{-1}})$, $\sup_{\Omega_2}|F_{H_0}|_{H_0}$, $\sup_{\Omega_2}|\theta|_{H_0}^{2}$, $\sup_{\Omega_2}|\sqrt{-1}\Lambda_\omega F_{H_0}|_{H_0}$,
$\sup_{\Omega_2}|D_{H_0}(\Lambda_{\omega }F_{H_0})|_{H_0}^2$, $\sup_{\Omega_2}|\nabla_{H_0}\theta|_{H_0}$, $d^{-2}$ and the geometry of $(\Omega_2, \omega)$. By choosing
$\tilde{B}=e^{C_{11}}(C_{12}+1)$,
then
\begin{equation} \label{noncompactc1eq4}
\sqrt{-1}\Lambda_\omega\partial\overline{\partial}\eta_i\geq \psi^2_2|T_i|^2_{H_i}-C_{13}
\end{equation}
on $M$. Let $\eta_i(P_0)=\sup_{M}\eta_i$. According to the definition of $\psi_{\alpha}$ and the uniform bounded on $\sup_{\Omega_2}\tr{h_i}$ , we may assume that $P_0\in \Omega_1$.
The inequality \eqref{noncompactc1eq4} and the maximum principle yield
\begin{equation*}
|T_i|^2_{H_i}(P_0)\leq C_{13},
\end{equation*}
and then there exists a positive constant $C_{14}$ depending only on $\sup_{\Omega_2}(\tr{h_i}+\tr{h_i^{-1}})$, $\sup_{\Omega_2}|F_{H_0}|_{H_0}$, $\sup_{\Omega_2}|\theta|_{H_0}^{2}$, $\sup_{\Omega_2}|\sqrt{-1}\Lambda_\omega F_{H_0}|_{H_0}$,
$\sup_{\Omega_2}|D_{H_0}(\Lambda_{\omega }F_{H_0})|_{H_0}^2$, $\sup_{\Omega_2}|\nabla_{H_0}\theta|_{H_0}$, $d^{-2}$ and the geometry of $(\Omega_2, \omega)$, such that
\begin{equation}
\sup_{\Omega }|T_i|^2_{H_0} \leq C_{14}.
\end{equation}
for all $i\geq i_{0}$. This concludes the proof of uniform local $C^1$-estimate of $H_{i}$.

We use Mochizuki's argument  (\cite{Mo3}) to give a uniform $L^2$ bounded on $|\partial_{H_0}h_i|_{H_{0}}=|\overline{\partial}_{E}h_{i}|_{H_{0}}$. Applying (\ref{key01}) and the Stokes formula, one can deduce
\begin{equation}\label{c1004}
\begin{split}
\int_{M_{i}} \sqrt{-1}\Lambda_\omega \overline{\partial}\partial \tr{h_i}\frac{\omega^{n}}{n!} &=\frac{1}{2}\int_{M_{i}} \sqrt{-1} (\overline{\partial}\partial -\partial \overline{\partial})\tr{h_i}\wedge\frac{\omega^{n-1}}{(n-1)!}\\
&=\frac{1}{2}\int_{M_{i}} \sqrt{-1} d(\partial -\overline{\partial})\tr{h_i}\wedge\frac{\omega^{n-1}}{(n-1)!}\\
&=-\frac{1}{2}\int_{M_{i}} d(\ast d\tr h_i)=-\frac{1}{2}\int_{M_{i}}div (\nabla \tr h_i)\frac{\omega^{n}}{n!} \\
&=-\frac{1}{2}\int_{\partial M_{i}}\frac{\partial \tr{h_i}}{\partial \nu_{i}} dv_{\partial M_i}\geq 0,
\end{split}
\end{equation}
where $\nu_{i}$ is the outer normal vector field at $\partial M_i$. This together with (\ref{c1003}) and the uniform bounded on $\tr{h_{i}}+\tr{h_{i}^{-1}}$ implies that there exists a uniform positive constant $C_{15}$ such that
\begin{equation}\label{l201}
\int_{M_{i}} (|  \overline{ \partial}_{E}h_i|_{H_{0}}^2 +|[\theta , h_i]|_{H_0}^{2})\frac{\omega^{n}}{n!}\leq C_{15}.
\end{equation}

Since we have obtained a uniform $C^{0}$-estimate and uniform local $C^{1}$-estimate on $h_i$, by the equation (\ref{c1002}) and the standard elliptic estimates, we can derive uniform local higher order estimates on $h_{i}$. So by choosing a subsequence (which we also denote by $H_i$), we know that $H_i$ converges to a Hermitian metric $H$ on whole $M$ in $C_{loc}^{\infty}$-topology as $i\rightarrow \infty$, and $H$ satisfies the Hermitian-Einstein equation (\ref{HEE}). Furthermore, (\ref{l201}) implies:  $\int_{M} (|  \overline{ \partial}_{E}h|_{H_{0}}^2 +|[\theta , h]|_{H_0}^{2})\frac{\omega^{n}}{n!}\leq C_{15}$. This complete the proof of Theorem \ref{theorem1}.
\vskip 1 true cm

\section{Hermitian-Einstein metrics}
In this section,  we follow Mochizuki's arguments to give a sufficient condition for the uniqueness of Hermitian-Einstein metric in Theorem \ref{theorem1}.


\begin{prop}(\cite{Mo3}, Proposition 2.4)\label{uniq}
Suppose $H_1$ and $H_2$ are two Hermitian-Einstein metrics of $(E,\overline{\partial}_E,\theta)$. Assume (1) $H_1$ and $H_2$ are mutually bounded, (2) $\sqrt{-1}\Lambda_\omega F_{H_1, \theta}=\sqrt{-1}\Lambda_\omega F_{H_2, \theta}$.
Then there exist a holomorphic decomposition $(E, \theta)=\bigoplus_{i=1}^m (E_i, \theta|_{E_i})$ and a tuple $(c_1, \ldots, c_m)\in \mathbb{R}_{>0}^m$ such that 
\begin{enumerate}
\item The decomposition $E=\bigoplus_{i=1}^m E_i$ is orthogonal with respect to both $H_i(i=1, 2)$,
\item $H_1|_{E_i}=c_i H_2|_{E_i}$.
\end{enumerate}
\end{prop}

\medskip

{\bf Proof} Let $h=H_1^{-1}H_2$. A direct computation gives us that
\begin{equation}
\sqrt{-1}\Lambda_\omega \partial\overline{\partial} \tr{h}=|  h^{-\frac{1}{2}}\cdot \partial_{H_1}h|_{H_{1}}^2 +|[\theta , h]\cdot h^{-\frac{1}{2}}|_{H_1}^{2}\geq 0.
\end{equation}
From Assumption 1, it can be seen that $\sqrt{-1}\Lambda_\omega \partial \overline{\partial} \tr{h}=0$. This means $|h^{-\frac{1}{2}}\cdot \partial_{H_1}h|_{H_{1}}=0$ and $|[\theta , h]\cdot h^{-\frac{1}{2}}|_{H_1}=0$. So $\partial_{H_1}h=0$ and $[\theta , h]=0$. Obviously the fact that $h$ is self-adjoint with respect to $H_i(i=1, 2)$ implies $\overline{\partial} h=0$ and $\partial_{H_2} h=0$. Then it follows that the eigenvalues of $h$ are constant. Let $E=\bigoplus_{i=1}^m E_i$ denote the eigen decomposition of $h$, which is the one we desired.

\begin{prop}
Let $(M, \omega )$ be a non-compact Gauduchon manifold satisfying the assumption 1 and $|d\omega^{n-1}|_{\omega} \in L^2(M)$. Suppose that there is a positive exhaustion function $\phi_1: M\rightarrow \mathbb{R}$ such that $|d\log\phi_1|_{\omega} \in L^2(M)$. Let $(E,\overline{\partial}_E,\theta)$ be an analytically $H_0$-stable bundle on $M$. Assume that $H_i(i=1,2)$ are Hermitian-Einstein metrics such that (1)$\det(H_i)=\det(H_0)$, (2) $H_i$ and $H_0$ are mutually bounded. Then $H_1=H_2$.
\end{prop}

{\bf Proof} Assume that $H_1$ satisfies Theorem \ref{theorem1} and let $h_1=H_0^{-1}H_1$.

Firstly Proposition \ref{uniq} gives us the decomposition $(E, \theta)=\bigoplus_{j=1}^m (E_j, \theta|_{E_j})$ such that (1) the decomposition is orthogonal with respect to $H_i (i=1, 2)$, (2) $H_1|_{E_j}=c_j H_2|_{E_j}$ for some $c_j > 0$. Let $\pi_j$ denote the projection onto $E_j$ with respect to the decomposition and $\pi_{j}^{\ast H_0}$ denote the adjoint of $\pi_j$ with respect to $H_0$. Because $H_0$ and $H_i (i=1, 2)$ are mutually bounded, one can immediately knows that $\pi_j$ are bounded with respect to $H_0$.

Set $\partial_{H_0, \theta}=\partial_{H_0}+\theta^{\ast H_0}$. Noting that $\theta^{\ast H_1}=h_1^{-1}\theta^{\ast H_0}h_1$ and $\partial_{H_1}-\partial_{H_0}=h_1^{-1}\partial_{H_0}h_1$, we deduce
\begin{equation}
\begin{split}\label{A}
\partial_{H_0}\pi_j=&\partial_{H_1}\pi_j-[h_1^{-1}\partial_{H_0}h_1, \pi_j]\\
=&\partial_{H_1}\pi_j-[h_1^{-1}(\partial_{H_0, \theta}h_1-[\theta^{\ast H_0}, h_1]), \pi_j]\\
=&\partial_{H_1}\pi_j-[h_1^{-1}\partial_{H_0, \theta}h_1, \pi_j]+[h_1^{-1}\theta^{\ast H_0}h_1-\theta^{\ast H_0}, \pi_j]\\
=&\partial_{H_1}\pi_j-[h_1^{-1}\partial_{H_0, \theta}h_1, \pi_j]+[\theta^{\ast H_1}, \pi_j]-[\theta^{\ast H_0}, \pi_j].\\
\end{split}
\end{equation}

%

By Mochizuki's arguments (\cite{Mo3}), we also consider the Hermitian metric $H_3=\bigoplus_{j=1}^m H_0|_{E_j}$. Then $H_3$ and $H_0$ are mutually bounded.

Set $h_3=H_0^{-1}H_3$. It can be expressed as that
\begin{equation}
h_3=\sum_{j=1}^{m}\pi_{j}^{\ast H_0}\circ \pi_j.
\end{equation}
And there holds that $\overline{\partial}_E h_3 \in L^2(H_0)$ and $\partial_{H_0}h_3\in L^2(H_0)$.

Denote $\overline{\partial}_{E, \theta}=\overline{\partial}_{E}+\theta$. Then we obtain
\begin{equation}\label{AAA}
\begin{split}
&\Lambda_{\omega}\overline{\partial}_{E, \theta}(h_3^{-1}\partial_{H_0, \theta}h_3)\\
=&\Lambda_\omega\overline{\partial}_E(h_3^{-1}(\partial_{H_0}h_3+[\theta^{\ast H_0}, h_3]))+\Lambda_\omega[\theta, h_3^{-1}(\partial_{H_0}h_3+[\theta^{\ast H_0}, h_3])]\\
=&\Lambda_\omega\overline{\partial}_E(h_3^{-1}\partial_{H_0}h_3)+ \Lambda_\omega\overline{\partial}_E(h_3^{-1}\theta^{\ast H_0}h_3-\theta^{\ast H_0})+\Lambda_\omega[\theta, h_3^{-1}\partial_{H_0}h_3]\\
&+\Lambda_\omega[\theta, h_3^{-1}\theta^{\ast H_0}h_3-\theta^{\ast H_0}]\\
=&\Lambda_\omega\overline{\partial}_E(h_3^{-1}\partial_{H_0}h_3)+\Lambda_\omega\overline{\partial}_E(\theta^{\ast H_3}-\theta^{\ast H_0})+\Lambda_\omega[\theta, h_3^{-1}\partial_{H_0}h_3]\\
&+\Lambda_\omega[\theta, \theta^{\ast H_3}-\theta^{\ast H_0}]\\
=&\Lambda_\omega\overline{\partial}_E(h_3^{-1}\partial_{H_0}h_3)+\Lambda_\omega[\theta, \theta^{\ast H_3}-\theta^{\ast H_0}]\\
=&\Lambda_\omega F_{H_3}-\Lambda_\omega F_{H_0}+\Lambda_\omega[\theta, \theta^{\ast H_3}]-\Lambda_\omega[\theta, \theta^{\ast H_0}].
\end{split}
\end{equation}

According to the holomorphic decomposition $(E, \theta)=\bigoplus_{i=1}^m (E_i, \theta|_{E_i})$ with respect to $H_i (i=1, 2)$, we have $\partial_{H_1, \theta}\pi_j=\partial_{H_1}\pi_j+ [\theta^{\ast H_1}, \pi_j]=0$ and $\overline{\partial}_{E, \theta}\pi_j=0$. Together with \eqref{A}, this means

\begin{equation}
\begin{split}
\partial_{H_0, \theta}\pi_j=&\partial_{H_0}\pi_j+[\theta^{\ast H_0}, \pi_j]\\
=&-[h_1^{-1}\partial_{H_0, \theta}h_1, \pi_j]-[\theta^{\ast H_0}, \pi_j]+[\theta^{\ast H_0}, \pi_j]\\
=&-[h_1^{-1}\partial_{H_0, \theta}h_1, \pi_j].
\end{split}
\end{equation}

Since $H_1$ is an Hermitian-Einstein metric on $(E, \theta)$, from Theorem \ref{theorem1}, we can get that $\partial_{H_0}h_1$, $\partial_{H_0, \theta}h_1$, $\partial_{H_0}\pi_j $ and $\partial_{H_0, \theta}\pi_j$ are in $L^2(H_0)$. Then it follows that $[\theta^{\ast H_0}, \pi_j] \in L^2(H_0)$ and $[\theta, \pi_j^{\ast H_0}] \in L^2(H_0)$. So immediately we know $\overline{\partial}_{E, \theta}\pi_j^{\ast H_0} \in L^2(H_0)$.

One can easily check that
\begin{equation}
\begin{split}
[\theta, h_3]=&\sum_{j=1}^{m}(\theta\circ\pi_{j}^{\ast H_0}\circ \pi_j-\pi_{j}^{\ast H_0}\circ \pi_j\circ\theta)\\
=&\sum_{j=1}^{m}(\theta\circ\pi_{j}^{\ast H_0}\circ \pi_j-\pi_j^{\ast H_0}\circ\theta\circ\pi_j+\pi_j^{\ast H_0}\circ\theta\circ\pi_j-\pi_{j}^{\ast H_0}\circ \pi_j\circ\theta)\\
=&\sum_{j=1}^{m}([\theta, \pi_j^{\ast H_0}]\circ\pi_j+\pi_j^{\ast H_0}\circ [\theta, \pi_j])\\
=&\sum_{j=1}^{m}[\theta, \pi_j^{\ast H_0}]\circ\pi_j.
\end{split}
\end{equation}
Then it can be seen that $[\theta, h_3] \in L^{2}(H_0)$ and $[\theta^{\ast H_0}, h_3] \in L^{2}(H_0)$. This implies that $\overline{\partial}_{E, \theta}h_3$ and $\partial_{H_0, \theta}h_3$ are square integrable with respect to $H_0$.

A direct computation gives us that
\begin{equation}
\overline{\partial}_{E, \theta}(h_3^{-1}\partial_{H_0, \theta}h_3)=-h_3^{-1}(\overline{\partial}_{E, \theta} h_3)h_3^{-1}\partial_{H_0, \theta}h_3+h_3^{-1}\overline{\partial}_{E, \theta}\partial_{H_0, \theta}h_3
\end{equation}
and
\begin{equation}
\overline{\partial}_{E, \theta}\partial_{H_0, \theta}h_3=\sum\overline{\partial}_{E, \theta}\pi_j^{\ast H_0}\circ \partial_{H_0, \theta}\pi_j+\sum\pi_{j}^{\ast H_0}\circ [F_{H_0, \theta}, \pi_{j}].
\end{equation}
Due to the assumption $|\Lambda_\omega F_{H_0, \theta}|_{H_0}\leq \hat{B}\varphi$, $|\Lambda_\omega F_{H_0, \theta}|_{H_0}$ is $L^1$. Then we obtain the following lemma.
\begin{lem}
$\Lambda_\omega \tr\overline{\partial}_{E, \theta}(h_3^{-1}\partial_{H_0, \theta}h_3)$ is $L^1$.
\end{lem}

Furthermore, we can derive
\begin{lem}\label{lemA}
$\int_M \tr(\overline{\partial}_{E, \theta}(h_3^{-1}\partial_{H_0, \theta}h_3))\wedge\omega^{n-1}=0$.
\end{lem}
{\bf Proof} Set $\chi_N:=\rho(N^{-1}\phi_1)$, where $\rho$ is a nonnegative $C^{\infty}$-function such that $\rho(t)=0$ if $t\geq 2$ and $\rho(t)=1$ if $t\leq 1$.. Because $\tr(\overline{\partial}_{E, \theta}(h_3^{-1}\partial_{H_0, \theta}h_3))$ is $L^1$, we just need to show
\begin{equation}\label{aaa}
\lim_{N\rightarrow \infty}\int_M \chi_N\cdot \tr(\overline{\partial}_{E, \theta}(h_3^{-1}\partial_{H_0, \theta}h_3))\wedge\omega^{n-1}=0.
\end{equation}
Computing directly yields that
\begin{equation}
\begin{split}
&\int_M \chi_N\cdot \tr(\overline{\partial}_{E, \theta}(h_3^{-1}\partial_{H_0, \theta}h_3))\wedge\omega^{n-1}\\
=&\int_M \chi_N\cdot \overline{\partial}(\tr(h_3^{-1}\partial_{H_0, \theta}h_3))\wedge\omega^{n-1}\\
=&\int_M\chi_N\cdot \tr(h_3^{-1}\partial_{H_0, \theta}h_3)\wedge\overline{\partial}\omega^{n-1}-\int_M\overline{\partial}\chi_N\cdot \tr(h_3^{-1}\partial_{H_0, \theta}h_3)\wedge\omega^{n-1}\\
=&\int_M\chi_N\cdot \partial(\log\det h_3)\wedge\overline{\partial}\omega^{n-1}-\int_{\{N\leq \phi_1 \leq 2N\}} \rho^{\prime}\big(\frac{\phi_1}{N}\big)\frac{\overline{\partial}\phi_1}{N}\cdot \tr(h_3^{-1}\partial_{H_0, \theta}h_3)\wedge\omega^{n-1}\\
=&-\int_M\partial{\chi_N}\cdot \log\det h_3\wedge\overline{\partial}\omega^{n-1}-\int_M\chi_N\cdot \log\det h_3\cdot\partial\overline{\partial}\omega^{n-1}\\
&-\int_{\{N\leq \phi_1 \leq 2N\}} \rho^{\prime}\big(\frac{\phi_1}{N}\big)\frac{\overline{\partial}\phi_1}{N}\cdot \tr(h_3^{-1}\partial_{H_0, \theta}h_3)\wedge\omega^{n-1}\\
=&-\int_{\{N\leq \phi_1 \leq 2N\}} \rho^{\prime}\big(\frac{\phi_1}{N}\big)\frac{\partial\phi_1}{N}\log\det h_3\wedge\overline{\partial}\omega^{n-1}\\
&-\int_{\{N\leq \phi_1 \leq 2N\}} \rho^{\prime}\big(\frac{\phi_1}{N}\big)\frac{\overline{\partial}\phi_1}{N}\cdot \tr(h_3^{-1}\partial_{H_0, \theta}h_3)\wedge\omega^{n-1}.
\end{split}
\end{equation}
Clearly $|d\omega^{n-1}|_{\omega} \in L^2(M)$ and $|d\log\phi_1|_{\omega} \in L^2(M)$ implies \eqref{aaa}.

Combining \eqref{AAA} and Lemma \ref{lemA}, we obtain
\begin{equation}
\int_M \Lambda_\omega \tr F_{H_0}=\int_M \Lambda_\omega \tr F_{H_3}=\sum_{i=1}^m\int_M\Lambda_\omega \tr F_{H_0|_{E_j}}.
\end{equation}
From $rank E= \sum_{j=1}^{m}rank E_j$, one can see that there exists $j_0$ such that $\mu(E, H_0)\leq \mu(E_{j_0}, H_0|_{E_{j_0}})$. This contradicts with the analytic stability of $(E, \overline{\partial}_E, H_0)$. So $H_1=H_2$.

\vskip 1 true cm


\bigskip
\bigskip

\noindent {\footnotesize {\it Chuanjing Zhang and Xi Zhang} \\
{School of Mathematical Sciences, University of Science and Technology of China}\\
{Anhui 230026, P.R. China}\\
{Email:chjzhang@mail.ustc.edu.cn;  mathzx@ustc.edu.cn}

\vskip 0.5 true cm

\end{document}